\documentclass[10pt]{amsart}
\usepackage{amssymb,latexsym,amscd}

\numberwithin{equation}{section}

\newtheorem{Thm}{Theorem}[section]

\newtheorem{Prop}[Thm]{Proposition}
\newtheorem{Lem}[Thm]{Lemma}
\theoremstyle{definition}
\newtheorem{Def}[Thm]{Definition}
\newtheorem{Ex}[Thm]{Example}
\newtheorem{Rmk}[Thm]{Remark}

\newtheorem{Conj}[Thm]{Conjecture}
\date{}

\newcommand{\LL}{{\mathcal{L}}} 
  
\newcommand{\II}{{\mathcal{I}}}   
\newcommand{\OO}{{\mathcal{O}}}
\newcommand{\PP}{{\mathbb{P}}}
\newcommand{\LP}{{\mathcal{L}}(-\sum_{i=1}^hm_iP_i)}
\newcommand{\ih}{i=1, \dots, h}
\newcommand{\bPP}{\tilde{\PP}}
\newcommand{\FF}{{\mathbb{F}}}

\begin{document}

\title{Special effect varieties and $(-1)$-curves}
\author{Cristiano Bocci}

\address{C. Bocci, Dipartimento di Matematica, Universit\`a di Milano,
Via Cesare Saldini 50, 20133 Milano, Italy}
\email{Cristiano.Bocci@unimi.it}
\thanks{This reseach was partially supported by GNSAGA of INdAM (Italy).}

\begin{abstract} 
Here we introduce the concept of special effect curve which permits to study, 
from a different point of view, special linear systems in $\PP^2$, i.e. 
linear system with general multiple base points whose effective dimension 
is strictly greater than the expected one. 
In particular we study two different kinds of special effect:
the $\alpha-${\it special effect} is defined
by requiring some numerical conditions,
while the definition of $h^1-${\it special effect}
concerns cohomology groups. 
We state two new conjectures for
the characterization of special linear systems
and we prove they are equivalent to the Segre and the
Harbourne-Hirschowitz ones.
\end{abstract}

\maketitle


\section{Introduction}
Let $X$ be a smooth, irreducible, complex projective variety of dimension $n$. 
Let $\LL$ be a complete linear system of divisors on $X$. 
Fix points $P_1,\dots, P_h$ on $X$ in general position and positive 
integers $m_1, \dots ,m_h$. We denote by $\LP$ the subsystem of $\LL$ given 
by all divisors having multiplicity at least $m_i$ at $P_i$, $i=1, \dots, h$. 
Since a point of multiplicity $m$ imposes 
$\left(\begin{smallmatrix} m+n-1 \\ n \end{smallmatrix}\right)$ conditions
we can define the {\bf virtual dimension} of the system $\LP$ as
$$\nu(\LP):=\mbox{virtdim}(\LP)=\dim(\LL)-\sum_{i=1}^h \binom{m_i+n-1}{n}.$$

This virtual dimension can be negative: in this case we expect that 
the system $\LP$ is empty. We can then define the 
{\bf expected dimension} of $\LP$ as
$$\epsilon(\LP):=\mbox{expdim}(\LP)=\mbox{max}\{\nu(\LP),-1\}.$$

The conditions imposed by the multiple points $m_iP_i$ can be 
dependent, so, in general we have 
\begin{center}
$\dim(\LP)\geq\epsilon(\LP)$
\end{center}
and we can state the following

\begin{Def}
A system $\LP$ is {\bf special} if 
\begin{center}
$\dim(\LP)>\epsilon(\LP)$, 
\end{center}
otherwise $\LP$ is said to be {\bf non-special}.
\end{Def}

By definition a system which is empty is non-special. 
For a non--empty system non-speciality means 
that the imposed conditions are independent. 

Since we expect that most systems are non-special, 
we can pose the following classification 
problem: {\it classify all special systems}.

The dimensionality problem is quite hard if we consider a 
general variety $X$, so we fix our attention on particular varieties and linear systems. 
As a first choice we can take $X=\PP^n$ and $\LL=\LL_{n,d}:=|{\mathcal{O}}_{\PP^n}(d)|$, 
the system of hypersurfaces of degree $d$ in $\PP^n$. In this case we have
\begin{center}
$\nu(\LL_{n,d}(-\sum_{i=1}^hm_iP_i))=\binom{d+n}{n} -1 -\sum_{i=1}^h \binom{m_i+n-1}{n}$.
\end{center}

Starting with the case $X=\PP^2$, we have some precise 
conjectures about the characterization of special linear systems and 
a rich series of results on the conjectures.
The main Conjectures are the following.
  
\begin{Conj}[(SC) B. Segre, 1961]\label{SC} If a linear system 
of plane curves with general multiple base points 
$\LL_{2,d}(-\sum_{i=1}^hm_iP_i)$ is special, then its general 
member is non-reduced, i.e. the linear system has, 
according to Bertini's theorem, some multiple fixed component.
\end{Conj}

\begin{Conj}[(HHC) Harbourne-Hirschowitz, 1989]\label{HHC} A linear system 
of plane curves with general multiple base points 
$\LL:=\LL_{2,d}(-\sum_{i=1}^hm_iP_i)$ is special if and only if is $(-1)-$special, 
i.e. its strict transform on the blow-up along the points $P_1, \dots,
P_h$ splits as 
$\tilde{\LL}=\sum_{i=1}^kN_iC_i+\tilde{{\mathcal{M}}}$
where the $C_i$, $i=1, \dots ,k$, are $(-1)-$curves such that 
$C_i\cdot \tilde{\LL}=-N_i<0$,  $\nu(\tilde{{\mathcal{M}}})\geq 0$ and 
there is at least one index $j$ such that $N_j>2$.
\end{Conj}

In \cite{CiMi3} C. Ciliberto and R. Miranda proved 
that the Harbourne--Hirschowitz and Segre Conjectures are equivalent.
Although the Harbourne--Hirschowitz Conjecture is still unproved, 
it is important to notice that, in more than a century of research, 
no special system has been discovered except $(-1)-$special systems. 
For an overview on these results the reader may consult \cite{Bocci},
 \cite{BoMi}, \cite{Ciliberto} and \cite{Miranda}.

When we pass to $\PP^n$, $n\geq 3$, very little is known 
about special linear systems.
One of the most important result
is the classification of the homogeneous special systems for double points:
\label{listaspeciali}
\begin{Thm}[Alexander--Hirschowitz, 1996]\label{P^nspecial}
The system ${\mathcal{L}}_{n,d}(2^h)$ is non-special unless:
$$\begin{array}{cccccc}
n  & any & 2 & 3 & 4 & 4 \cr
d  & 2   & 4 & 4 & 4 & 3 \cr 
h  & 2, \dots, n & 5 & 9 & 14 & 7 \cr
\end{array}$$
\end{Thm} 

Continuing with $\PP^n$, $n\geq3$
we can notice that there is not a precise conjecture.
Although the Segre Conjecture can be generalized
in every ambient variety using the statement concerning 
$H^1\not=0$ (see, for example, \cite{Bocci} or \cite{CiMi3})
there is nothing that characterizes the special systems
from a geometric point of view as, for example,
in the case of $(-1)-$curves in $\PP^2$.

A worthy goal would be ``find a conjecture (C) in $\PP^n$,
[or in a generic variety $X$] such that, when we read (C) in $\PP^2$,
(C) is equivalent to the Segre (\ref{SC})
and Harbourne--Hirschowitz (\ref{HHC}) Conjectures''.

In Sections \ref{alphasec} and \ref{h1sec} we state two potential
candidates for the above-mentioned goal: the {\it Numerical Special Effect Conjecture}
and the {\it Cohomological Special Effect Conjecture}.
In fact, in these sections, we define the concepts of
``$\alpha$'' and ``$h^1$'' special effect curves which permit to 
introduce a different approach in the study of special linear systems
in $\PP^2$. Moreover, in Section \ref{equivalence} we prove that
these conjectures are equivalent to the Segre and the
Harbourne-Hirschowitz ones. 

In Section \ref{highexamples} we present some examples of special
effect varieties in $\PP^n$, $n \geq 3$. Due to its complexity, the generalization of the
``Numerical'' and ``Cohomological'' Conjectures
to the higher dimensional case is presented in \cite{higherSEV} where
we prove also that
these Conjectures hold for every special system listed in Theorem
\ref{P^nspecial}. 

Finally, in section \ref{surfaces} we show some results on special effect
varieties when the ambient variety is a Hirzebruch surface or a $K3$
surface.
\vskip0.5cm
The main ideas of this article were born during a pleasant stay in
Fort Collins. I'm very grateful to Professor Rick Miranda for the
guidance and support during my research at the Colorado State
University and to
all the people that I met there, especially the staff and all
Professors at the Department of Mathematics.

I also thank Professors Luca Chiantini and Ciro Ciliberto for discussing topics about special effect varieties and giving me several interesting suggestions.


\section{Preliminaries}

We collect some facts about linear systems that will
be useful in the next Sections.

Consider the blow-up $\pi:\bPP^n \to \PP^n$ at the points $P_1, \dots, P_h$ 
and let $E_i$, $\ih$ be the exceptional divisors corresponding 
to the blow-up of the points $P_i, \ih$. 
If we denote by $H$ the pull-back of a general hyperplane of $\PP^n$ via $\pi$, 
then we can write the strict transform of the system 
$\LL:=\LL_{n,d}(\sum_{i=1}^hm_iP_i)$ 
as $\tilde{\LL}=|dH-\sum_{i=1}^hm_iE_i|$. 
In the future, if confusion cannot arise, we will indicate both $\LL$ and $\tilde{\LL}$ by $\LL$.

It an easy application of the (generalized) Riemann-Roch theorem to observe that
\begin{equation}\label{chi}
\nu(\LL)=\chi(\tilde{\LL})-1.
\end{equation}

Consider now the case of $\PP^2$ and let $\LL:=\LL_{2,d}(\sum_{i=1}^hm_iP_i)$. By Riemann--Roch, remembering that $h^2(\bPP^2,\tilde{\LL})=0$, we obtain
\begin{equation}\label{legami}
\begin{split}
\dim(\LL)=&\dim(\tilde{\LL})=\frac{\tilde{\LL} \cdot(\tilde{\LL}-\tilde{K})}{2}+h^1(\bPP^2,\tilde{\LL})-h^2(\bPP^2,\tilde{\LL})=\\
=&\tilde{\LL}^2-g_\LL+1+h^1(\bPP^2,\tilde{\LL})=\nu(\LL)+h^1(\bPP^2,\tilde{\LL})
\end{split}
\end{equation}
where $g$ is the arithmetic genus $p_a$ of a curve in $\tilde{\LL}$ and $\tilde{K}$ is the canonical class on $\bPP^2$. 

Hence, by previous formula, we have 
\begin{equation}\label{h0h1}
\LL \text{ is non-special if and only if } h^0(\bPP^2,\tilde{\LL})\cdot h^1(\bPP^2,\tilde{\LL})=0.
\end{equation}

\begin{Rmk}\rm The reducible curve $C=\sum_{i=1}^kN_iC_i$ 
in Conjecture \ref{HHC} is called a $(-1)-${\bf configuration} on $\bPP^2$.
\end{Rmk}

Whenever not otherwise specified, we work over the field ${\mathbb{C}}$. 


\section{$\alpha-$special effect curves}\label{alphasec}

Let $P_1, \dots, P_h$ be points in $\PP^2$ in general position 
and fix positive integers $m_1, \dots, m_h$. 
Consider the system $\LL:=\LL_{2,d}(-\sum _{i=1}^hm_iP_i)$ of planar 
curves of degree $d$ passing through the points $P_i$ 
with multiplicity at least $m_i$.

\begin{Def}\label{aSEP} 
Let $\LL$ and $P_i, \dots P_h$ as above. 
An irreducible curve $Y$, of degree $e$, has the 
$\alpha-${\bf special effect property} for $\LL$ on $\PP^2$ 
if there exist non-negative integers $\alpha, c_{j_1}, \dots c_{j_s}$, 
with $\alpha e \leq d$ and $1 \leq \alpha \leq \mbox{min}\{\lceil \frac{m_{j_i}}{c_{j_i}}\rceil, i=1, \dots, s\}$, 
such that
\begin{itemize}
\item[(i)]  $Y$ contains the point $P_{j_i}$ with multiplicity at least
  $c_{j_i}$ for $j=1, \dots, s$, where $P_{j_i}\in \{ P_1, \dots, P_h \}$; 
\item[(ii)] $\nu(\LL-\alpha Y)> \nu(\LL)$.
\end{itemize}
Moreover we require that $\alpha$ is the maximum admissible value for
the $\alpha-$special effect property and, if $\beta > \alpha$ then
$\nu(\LL-\beta Y)<\nu(\LL-\alpha Y)$.
\end{Def}

In the following, we will mainly ask for a condition stronger than $(i)$:
\begin{itemize}
\item[(i*)] $\nu(|Y|)\geq 0$, 
\end{itemize}
where $|Y|$ represents the linear system
$|Y|=|eH-\sum_{i=1}^sc_{j_i}P_{j_i}|$. It is clear that condition $(i*)$ implies condition
$(i)$. 

Condition $(ii)$ is surely the most interesting. 
As a matter of fact it tells us that the number of conditions 
imposed on the system of curves of degree $d$ by imposing a multiple 
curve $\alpha Y$ and the points $P_{j_i}$ with multiplicity $m_{j_i}-\alpha c_{j_i}$ 
(such that the final multiplicity at the point $P_{j_i}$ is at least
$m_{j_i}$, $i=1, \dots s$) plus eventually the other multiple points
$m_tP_{t}$, $t \not\in \{ j_1, \dots, j_s \}$
is less than the number of conditions imposed to the same system $|dH|$ 
only imposing each $P_i$ with multiplicity at least $m_i$, $\ih$.  
This sounds like a crazy requirement because, in general, we expect that a 
positive dimensional variety imposes more conditions than a zero-dimensional variety. 
It is important to notice the similarity with the ``strange'' 
requirement in the case of $(-1)-$curves in \cite{Ciliberto}: 
we asked there for a curve $C$ whose double is not expected to exist !

\begin{Ex}\label{L12}\rm
Let $\LL$ be the system $\LL_{2,9}(-6P_1-6P_2-6P_3)$.
This system is special since $\nu(\LL)=-9$ but its effective dimension is $0$ 
since it contains $3Y$, with $Y=L_{12}+L_{13}+L_{23}$,
where $L_{ij}$ is the line through $P_i$ and $P_j$.
We claim that each of the lines $L_{ij}$
has the $3-$special effect property.
We prove this for $L_{12}$.
Obviously one has $\nu(|L_{12}|)\geq0$; indeed, it is a $(-1)$-curve.
Moreover $\LL-L_{12}$ is the system
$\LL':=\LL_{2,8}(-5P_1-5P_2-6P_3)$
and its virtual dimension is 
\[
\nu(\LL'):=\frac{8\cdot 11}{2}-2\frac{5\cdot6}{2}-\frac{6\cdot7}{2}=44-30-21=-7.
\]
Going further we can observe that
\[
\nu(\LL-2L_{12})=\nu(\LL-3L_{12})=-6
\]
while
\[
\nu(\LL-4L_{12})=-7.
\]
So the claim follows. 
\end{Ex}

\begin{Ex}\rm
Let $\LL:=\LL_{2,d}(-\sum_{i=1}^h m_iP_i)$
and consider a $(-1)-$curve $E$ such that $\LL\cdot E=-N<0$.
Thus $\LL=NE+{\mathcal{M}}$, where $E\cdot{\mathcal{M}}=0$.
Using Riemann-Roch it is easy to prove
$\nu(\LL-NE)=\nu(\LL)+\binom{N}{2}$ and $\nu(\LL-(N+1)E)=\nu(\LL-NE)-1$.
Hence $E$ has the $N-$special effect property if $N\geq 2$.
\end{Ex}

Going back to the definition of $\alpha-$special effect curves, 
we now see how the conditions $(i)-(ii)$ give some numerical information 
about the intersection $\LL\cdot Y$.  We will also work on the blow-up of 
$\PP^2$ at the points $P_1, \dots, P_h$ and, as in the case of $(-1)-$curves, 
we will consider the strict transform $\tilde{Y}$ of the
$\alpha-$special effect curve $Y$, but in general we will denote 
both $Y$ and $\tilde{Y}$ by $Y$. 

\begin{Lem}\label{LY}
Let $Y$ be an irreducible curve having the $\alpha-$special effect property 
for a system $\LL$. Then $\LL \cdot Y < \frac{(\alpha+1)}{2}Y^2$.
\end{Lem}

\begin{proof}
Let $\LL$ be the system $\LL_{2,d}(-\sum_{i=1}^hm_iP_i)$ 
and suppose $Y$ has degree $e$ and passes through $P_{j_i}$'s with multiplicity at least $c_{j_i}$.  
From conditions $(i*)$ and $(ii)$ of the $\alpha-$special effect
property we have respectively
\begin{alignat}{1}
e^2+3e\geq \sum_{i=1}^s(c_{j_i}^2+c_{j_i})  \mbox{ then } -3e+\sum_{i=1}^sc_{j_i}\leq e^2-\sum_{i=1}^sc_{j_i}^2 \label{1}\\
\frac{1}{2}\left(-de+ \sum_{i=1}^sm_{j_i}c_{j_i}+\alpha e^2-3e- \sum_{i=1}^s(\alpha c_{j_i}^2+c_{j_i})\right)>0\label{2}
\end{alignat}
Since $\LL\cdot Y:=\tilde{\LL}\cdot\tilde{Y}=de-\sum_{i=1}^sm_{j_i}c_{j_i}$, we obtain (by using (\ref{2}) and (\ref{1})): 
$$\LL\cdot Y = de-\sum_{i=1}^sm_{j_i}c_{j_i}<\frac{1}{2}\left( \alpha e^2-3e- \sum_{i=1}^s(\alpha c_{j_i}^2+c_{j_i})  \right) \leq \frac{(\alpha+1)}{2}(e^2- \sum_{i=1}^sc_{j_i}^2)$$
so that $\LL \cdot Y < \frac{(\alpha+1)}{2}Y^2$.
\end{proof}

By the previous lemma we can also obtain some informations about $Y^2$.

\begin{Lem}\label{inter}
Suppose $Y$ has the $\alpha-$special effect property for a system $\LL$. 
If $h^0(\LL-\alpha Y)\geq 1$ then $Y^2\leq -1$.
\end{Lem}

\begin{proof}
By Lemma \ref{LY} we have
\begin{equation}\label{splitting}
(\LL-\alpha Y)\cdot Y=\LL\cdot Y-\alpha Y^2<\frac{(1-\alpha)}{2}Y^2
\end{equation}
Consider first the case $\alpha=1$; then $Y$ splits from $\LL-Y$ and we can compute
$$(\LL-2Y)\cdot Y= \LL\cdot Y - 2Y^2=\LL\cdot Y -Y^2-Y^2<-Y^2$$
Hence, if $Y^2 \geq 0$ then $Y$ is a fixed component of $\LL-2Y$. But at this point we can 
iterate the procedure and we would obtain
$$(\LL-NY)\cdot Y= \LL\cdot Y -Y^2-(N-1)Y^2<-(N-1)Y^2$$
Thus if $Y^2\geq 0$, $Y$ appears with multiplicity $\infty$ in $\LL-Y$, but this is a 
contradiction, hence $Y^2\leq -1$.

Consider now the case $\alpha \geq 2$ in (\ref{splitting}). If $Y^2 \geq 0$, then $Y$ is a 
fixed component of $\LL-\alpha Y$. 
Moreover, for $N>\alpha$, we have
$$ (\LL-N Y)\cdot Y=\LL\cdot Y-N Y^2<\frac{(\alpha+1-2N)}{2}Y^2<0.$$
Thus we can conclude again that if $Y^2\geq 0$, then we obtain a contradiction. Hence $Y^2\leq -1$.
\end{proof}

\begin{Def}\label{aSEV} 
Let $\LL$ and $P_1, \dots, P_h$ as above. An irreducible curve $Y$, of degree $e$, is 
an $\alpha-${\bf special effect curve} for $\LL$ on $\PP^2$ if $Y$ has the 
$\alpha-$special effect property for $\LL$ and moreover $\nu(\LL-\alpha Y)\geq 0$.
\end{Def}

We recall that the existence of a $(-1)-$configuration $C=\sum_{i=1}^tN_iC_i$ such that $\LL:=\sum_{i=1}^tN_iC_i+{\mathcal{M}}$ leads us to the inequality
\begin{equation}\label{NC}
\dim(\LL)=\dim({\mathcal{M}})\geq \nu({\mathcal{M})}=\nu(\LL)+\sum_{i=1}^t\binom{N_i}{2}.
\end{equation}
which, under the assumption of $(-1)-$speciality of $\LL$,
i.e. $\nu({\mathcal{M}})\geq 0$ and $N_i\geq 2$ for at least one index
$i$, implies that $\LL$ is special. 
Observe that the existence of an $\alpha-$special effect curve $Y$ for a system $\LL$ forces the system itself to be special. 
In fact  we have the following chain of inequalities
$$\dim(\LL)\geq\dim(\LL-\alpha Y )\geq\nu(\LL-\alpha Y)>\nu(\LL)$$
and, together with condition $\nu(\LL-\alpha Y)\geq 0$, one has $\dim(\LL)>\epsilon(\LL)$.

\begin{Ex}\rm
Let $\LL:=\LL_{2,2}(-2P_1-2P_2)$ be the linear system of conics with two double points. 
Let $Y$ be a line through $P_1$ and $P_2$, i.e $Y=H-P_1-P_2$. Obviously condition $(i)$ is satisfied. 
Since 
\[
\nu(\LL-Y)=\nu(\LL-2Y)=0
\]
while $\nu(\LL)=-1$, one has that condition $(ii)$ is satisfied. 
From the positivity of $\nu(\LL-2Y)$ we conclude that the line through $P_1$ and $P_2$ is a $2-$special 
effect curve for $\LL$ and so $\LL$ is special.  
\end{Ex}

\begin{Ex}\label{e1e2}\rm
We want to show how the problem of the existence of an 
$\alpha-$special effect curve can turn into a pure combinatorial problem and 
its solution is more or less difficult according to the initial data. 

For example, we can look for an irreducible smooth $\alpha-$special
effect curve $Y$ of degree $e$ for a generic homogeneous system $\LL:=\LL_{2,d}(m^h)$. Moreover we
require that $Y$ passes through all points $P_1, \dots, P_h$.
The smoothness of $Y$ means $c_1=\dots=c_h=1$. 

The conditions for the existence of $Y$ are:  
\begin{itemize}
\item[(i)] $P_i \in Y \mbox{ for } i=1, \dots , h$,
\item[(ii)] $\nu(|(d-\alpha e)H- \sum_{i=1}^h(m-\alpha)P_i|) > \nu(|dH-\sum_{i=1}^n mP_i|)$.
\item[(iii)] $\nu(|(d-\alpha e)H-\sum_{i=1}^h(m-\alpha)P_i|)\geq 0$, 
\end{itemize}
with the extra conditions $1 \leq \alpha \leq m$ and $\alpha e \leq d$. 
Using Riemann--Roch we can write the previous conditions as
\begin{alignat}{2}
& \frac{e(e+3)}{2}\geq h \label{Yexists} \\ 
& \frac{(d-\alpha e)(d+\alpha e+3)}{2}-h\frac{(m-\alpha)(m-\alpha+1)}{2} > \frac{d(d+3)}{2}- h\frac{m(m+1)}{2} \label{sineq1}\\
& \frac{(d-\alpha e)(d+\alpha e+3)}{2}\geq h\frac{(m-\alpha)(m-\alpha+1)}{2} \label{dimension}.
\end{alignat}
In particular, if we expand condition (\ref{sineq1}), we obtain
\begin{equation}\label{sineq} 
-d\alpha e + \frac{1}{2}\alpha^2e^2-\frac{3}{2}\alpha e +hm\alpha - \frac{1}{2}h \alpha^2 + \frac{1}{2}h \alpha > 0.
\end{equation}

Observe that (\ref{dimension}) is increasing monotone in $d$ and for $d=\alpha e$, we have
$$0 - h \binom{m-\alpha+1}{2} \geq 0$$
which is satisfied only for $\alpha = m $. Then $d=me$.

We claim that $d\geq me$. The proof of this fact is a 
long and very tedious study of the equations
{\small{\[
m^2e^2-2met-2me^2\alpha+3me+t^2+2t\alpha e-3t+\alpha^2e^2-3\alpha e-hm^2+2hm\alpha-hm-h\alpha^2+h\alpha  \geq 0
\]}}
and
\[
-2me^2\alpha+2t\alpha e + \alpha^2 e^2-3\alpha e +2hm\alpha -h\alpha^2+h\alpha>0\]
given by (\ref{dimension}) and (\ref{sineq}) in which we substitute $d=me-t$, with $t>0$.
Anyway, the previous equations together with $e^2+3e\geq 2h$ are 
verified only if at least one between $m,e,t,h$ and $\alpha$ is equal to zero, 
but this is not acceptable for our purposes (we can check it by a computer algebra system, e.g. Maple). 

Now we show that $e<3$: using (\ref{sineq}) we compute
$$d\alpha < \frac{1}{2}\alpha^2e-\frac{3}{2}\alpha+\frac{h\alpha}{e}(m-\frac{1}{2}\alpha+\frac{1}{2})$$
and from $d\geq me$ and (\ref{Yexists}) we obtain
$$m\alpha e \leq d\alpha < \frac{1}{2}\alpha^2e-\frac{3}{2}\alpha+\alpha(m-\frac{1}{2}\alpha+\frac{1}{2})\frac{(e+3)}{2}.$$
Then, 
$$m e < \frac{1}{2}\alpha e-\frac{3}{2}+\frac{1}{2}me-\frac{1}{4}\alpha e+\frac{1}{4} e+3m-\frac{3}{4}\alpha +\frac{3}{4}$$
and simplifying, we obtain
$$e(\frac{1}{2}m-\frac{1}{4}\alpha -\frac{1}{4})<3(\frac{1}{2}-\frac{1}{4}\alpha -\frac{1}{4})$$
that is $e<3$.

If we analyze the cases $e=1$ and $e=2$, we see that the only possibilities are 
\begin{center}
\begin{itemize}
\item $e=1$, $h=2$, $m\leq d <2m-\frac{1}{2}-\frac{1}{2}\alpha$
\item $e=2$, $h=5$, $2m\leq d <\frac{5}{2}m-\frac{1}{4}-\frac{1}{4}\alpha$
\end{itemize}
\end{center}
If we substitue $\alpha:=\LL \cdot Y = de-hm$ we obtain
\begin{center}
\begin{itemize}
\item $e=1$, $h=2$, $m\leq d <2m-2$
\item $e=2$, $h=5$, $2m\leq d <\frac{5m-2}{2}$
\end{itemize}
\end{center}

Then we conclude that the systems
$$\begin{array}{ll}
\LL_{2,d}(m^2) & m\leq d <2m-2 \\
\LL_{2,2d}(m^5) & 2m\leq d <\frac{5m-2}{2}
\end{array}$$
are special. The careful reader can observe that these families of 
special systems are exactly the first two cases in the 
classification of the homogeneous $(-1)-$special systems described in Theorem 2.4 in \cite{CiMi2}.
\end{Ex}

\begin{Rmk}\label{12}\rm Let $\LL$ be again the system $\LL_{2,9}(-6P_1-6P_2-6P_3)$. 
As already saw in Example \ref{L12} we know that each of the lines $L_{ij}$ has the 
$3-$special effect property for $\LL$. As we can see, a single line is not a $3-$special 
effect curve for $\LL$, since $\nu(\LL-3L_{ij})<0$.
\end{Rmk}

The previous Remark shows that $\alpha-$special effect curves are 
not sufficient to describe all known special systems. Hovewer it is clear, now, 
in which way we proceed. If $Y$ has the $\alpha-$special effect property for a 
system $\LL$ and $\nu(\LL-\alpha Y)<0$, we substitute the system $\LL$ 
with $\LL-\alpha Y$ and we investigate this new system.

\begin{Def}\label{aSEC}
Let $\LL$ be a system as above.
Fix a sequence of (not necessarily distinct)
irreducible curves $Y_1, \dots Y_t$,
Suppose further that
\begin{itemize}
\item[(1)] $Y_j$ has the $\alpha_j-$special effect property 
for $\LL-\sum_{i=1}^{j-1}\alpha_iY_i$, for $j=1,\dots, t$,
\item[(2)] $\nu(\LL-\sum_{i=1}^t\alpha_i Y_i)\geq 0$.
\end{itemize}
Then we call both $X:=\sum_{i=1}^t\alpha_i Y_i$ and $\{ Y_1, \dots,
Y_t\}$ an $(\alpha_1, \dots, \alpha_t)-${\bf special effect configuration for
$\LL$}.
\end{Def}

\begin{Ex}\rm
Consider again the system $\LL:=\LL_{2,9}(-6P_1-6P_2-6P_3)$. 
We prove now that $X=3L_{12}+3L_{13}+3L_{23}$ is a $(3,3,3)-$special effect configuration. 
Recall that $\nu(\LL)=-9$. In Example \ref{L12} we proved that $L_{12}$ 
has the $3-$special effect property for $\LL$. We can go ahead and
check if $L_{13}$ 
has the $3-$special effect property for $\LL-3L_{12}$. We obtain:
\vskip0.2cm
$\begin{array}{l}
\nu(\LL-3L_{12}-L_{13})=\nu(|5H-2P_1-3P_2-5P_3|)=-4 \\
\nu(\LL-3L_{12}-2L_{13})=\nu(|4H-P_1-3P_2-4P_3|)=-3 \\
\nu(\LL-3L_{12}-3L_{13})=\nu(|3H-3P_2-3P_3|)=-3 
\end{array}$
\vskip0.2cm
\noindent Finally we check if $L_{23}$ has the $3-$special effect property for $\LL-3L_{12}-3L_{13}$:
\vskip0.2cm
$\begin{array}{l}
\nu(\LL-3L_{12}-3L_{13}-L_{23})=\nu(|2H-2P_2-2P_3|)=-1 \\
\nu(\LL-3L_{12}-3L_{13}-2L_{23})=\nu(|H-P_2-P_3|)=0.\\
\nu(\LL-3L_{12}-3L_{13}-3L_{23})=0.\\
\end{array}$
\vskip0.2cm
\noindent Thus  $X$ is a $(3,3,3)-$special effect configuration for $\LL_{2,9}(-6P_1-6P_2-6P_3)$.
\end{Ex}

As in the case of $\alpha-$special effect curves also a 
special effect configuration $X$ forces a system to be special. 
In fact, one has again 
$$\dim(\LL)\geq\dim(\LL-X )\geq\nu(\LL-X)>\nu(\LL)$$
and, together with condition $(2)$ in Definition \ref{aSEC}, one has $\dim(\LL)>\epsilon(\LL)$.

These facts permit us to define a particular kind of speciality. 

\begin{Def}
A special system arising from the existence of an 
$\alpha-$special effect curve (or an $(\alpha_1, \dots, \alpha_r)-$special 
effect configuration) is called {\bf Numerically Special}.
\end{Def}

Finally, we can state the following

\begin{Conj}[(NSEC) ``Numerical Special Effect'' Conjecture]\label{NSEC}
A linear system of plane curves $\LL_{2,d}(-\sum_{i=1}^hm_iP_i)$ 
with general multiple base points is special if and only if it is numerically special.
\end{Conj}


\section{$h^1-$Special effect curves}\label{h1sec}
The second class of curves we introduce are defined via some 
particular conditions on certain cohomology groups. The original idea 
for these curves comes from a detailed analysis of the base locus in the 
special systems listed in Theorem \ref{P^nspecial}, that is, linear systems 
with imposed double points in $\PP^n$, $n\geq 2$. In fact, as shoved
in \cite{higherSEV}, this kind of speciality can be more easily generalized to 
higher dimensions than {\it numerical} one.
 
\begin{Def}\label{h1SEV} Let $\LL:=\LL_{2,d}(-\sum_{i=1}^hm_iP_i)$ 
be a linear system of plane curves with general multiple base points.
An irreducible curve $Y\subset \PP^2$, with $\OO_{\PP^2}(Y)\not\cong \LL$, 
is an $h^1-${\bf special effect curve} for the system  $\LL$  if the following 
conditions are satisfied:
\begin{itemize}
\item[(a)] $h^0(\LL_{|Y})=0$;
\item[(b)] $h^0(\LL-Y)>0$;
\item[(c)] $h^1(\LL_{|Y})>0$.
\end{itemize}
\end{Def}

\begin{Rmk} \rm The condition $(c)$ will be slightly different in the
  definition in the higher dimension case where we ask for $h^1(\LL_{|Y})>h^2(\LL-Y)$. 
Instead, in the planar case, we can just ask for $h^1(\LL_{|Y})>0$ because $h^2(\LL-Y)=0$. 
In fact, by definition of $Y$ and condition $(b)$ we can suppose 
$\LL-Y=|aH-\sum_{i=1}^hs_iP_i|$, with $a, s_1, \dots, s_h$ positive integers. Define now $Z$ 
as the union of the fat points $s_iP_i$, then we have the following exact sequence
$$0 \to \II_Z\otimes \OO_{\PP^2}(a) \to  \OO_{\PP^2}(a) \to \OO_Z \to 0$$
When we consider the cohomology groups, we have
$$
\begin{array}{ccc}\dots \to h^1(\OO_Z) \to & h^2(\II_Z\otimes \OO_{\PP^2}(a)) &  \to h^2(\OO_{\PP^2}(a)) \to \dots\\
& \parallel \\
& h^2(\LL-Y) 
\end{array}$$
Since $Z$ is a zero-dimensional scheme one has $h^i(\OO_Z)=0$ for $i\geq 1$. 
Moreover, by Serre duality, $h^2(\OO_{\PP^2}(a))=h^0(\OO_{\PP^2}(-3-a))=0$. Thus $h^2(\LL-Y)=0$.
\end{Rmk}

\begin{Ex}\label{h1222}\rm
Let $\LL:=\LL_{2,2}(-2P_1-2P_2)$ be the linear system of conics 
with two double points. Let $Y$ be a line through $P_1$ and $P_2$, i.e $Y=H-P_1-P_2$. 
Since $\LL\cdot Y=-2$ the restricted system $\LL_{|Y}$ 
has no effective divisors and $h^0(\LL_{|Y})$ is empty. 
By Riemann--Roch we easily compute $h^1(\LL_{|Y})=g_Y-1-\deg(\LL_{|Y})=1>0$. 
Finally $\LL-Y$ is $|H-P_1-P_2|$, so that $h^0(\LL-Y)=1$. 
Hence the line $Y$ through $P_1$ and $P_2$ is an $h^1-$special effect curve for $\LL$.
\end{Ex}

Let $\LL:=\LL_{2,d}(-\sum_{i=1}^h m_iP_i)$ and consider,
on the blow-up of $\PP^2$ at the points $P_i$'s,
the exact sequence
\[
0 \to \LL-Y \to \LL \to \LL_{|Y} \to 0
\]
which gives the following long exact sequence in cohomology:
\[
0 \to H^0(\LL-Y) \to H^0(\LL) \to H^0(\LL_{|Y})
\to  H^1(\LL-Y) \to H^1(\LL) \to H^1(\LL_{|Y}) \to 0.
\]
Conditions $(a)$ and $(b)$
assure us that $H^0(\LL)\not= 0$,
while condition $(c)$ implies $H^1(\LL)\not= 0$.
Thus the existence of such $Y$
forces the system $\LL$ to have $h^0(\LL)\cdot h^1(\LL)\not= 0$
so that, by (\ref{h0h1}), $\LL$ is special.
Again, we can give a particular name to this kind of system:

\begin{Def}
A special system arising from the existence of an $h^1-$special 
effect curve is called {\bf Cohomologically Special}.
\end{Def}

And again we can state a conjecture:

\begin{Conj}[(CSEC) ``Cohomological Special Effect'' Conjecture]\label{CSEC}
A linear system of plane curves  $\LL:=\LL_{2,d}(-\sum_{i=1}^hm_iP_i)$ 
with general multiple base points is special if and only if it is cohomologically special.
\end{Conj}

\begin{Lem}
Suppose that the Cohomological Special Effect Conjecture holds. 
Let $C\subset \PP^2$ be an irreducible curve passing through the 
general points $P_1,\dots, P_h$ with multiplicity at least $m_1, \dots, m_h$. 
Then $\tilde{C}^2\geq g_{\tilde{C}}-1$.
\end{Lem}

\begin{proof}
Suppose $\nu(|C|)<0$, then the system $|C|$ is special. 
Thus there is an $h^1-$special effect curve $Y$ for $|C|$ and $Y$ is a fixed part of $C$. 
This is a contradiction since $C$ is irreducible. 
Hence $\nu(|C|)\geq 0$ and, by formula (\ref{legami}), one has $\tilde{C}^2\geq g_{\tilde{C}}-1$. 
\end{proof}


\section{The four conjectures}\label{equivalence}
In the previous sections we introduced two new conjectures 
for the characterization of special linear systems in the planar case.
At this point it is natural to ask if these conjectures are 
equivalent to the Segre and Harbourne--Hirschowitz ones. 
The answer is given in the following
\begin{Thm}
Conjectures {\rm (SC)}, {\rm (HHC)}, {\rm (NSEC)} and {\rm (CSEC)} are equivalent.
\end{Thm}

\begin{proof}
First of all, we recall that the equivalence between (SC) and (HHC) 
is proved in \cite{CiMi3}. Then we just need to prove the following implications:
$$\begin{matrix}
\mbox{(HHC)} \Rightarrow \mbox{(NSEC)}  \Rightarrow \mbox{(SC)} \\
\mbox{(HHC)} \Rightarrow \mbox{(CSEC)}  \Rightarrow \mbox{(SC)} \\
\end{matrix}
$$
\vskip0.3cm
\noindent
{\bf[(HHC) $\Rightarrow$ (NSEC)]} 
Suppose that the Harbourne--Hirschowitz Conjecture holds. 
Let $\LL$ be a special system, then it splits as $\LL=\sum_{i=1}^tN_iC_i+ {\mathcal{M}}$, 
where $\nu({\mathcal{M}})\geq 0$ and there is at least one index $j$ such that $N_j>1$.
After a permutation in the indexes we can suppose that $N_i>1$ for $i=1,\dots, s$, $s\leq t$. 
Thus we can write the $(-1)-$configuration $C=\sum_{i=1}^tN_iC_i$ appearing in $\LL$ as 
$$C=\sum_{i=1}^sN_iC_i+\sum_{i=s+1}^tC_i.$$
At this point it is enough to show that 
$\sum_{i=1}^sN_iC_i$ is an $(N_1, \dots, N_s)-$special effect configuration for $\LL$.
By formula (\ref{NC}) each $(-1)-$curve $C_j$ with $N_j>1$ increases the virtual 
dimension of the residual system by
$$\nu(\LL-\sum_{i=1}^jN_iC_i)=\nu(\LL-\sum_{i=1}^{j-1}N_iC_i)+\binom{N_i}{2}$$
Thus, if $N_j>1$ then $C_i$ has the $N_j-$special 
effect property for $\LL-\sum_{i=1}^{j-1}N_iC_i$.
Finally, we can observe that ${\mathcal{M}}=\LL-\sum_{i=1}^tN_iC_i$. 
By hypothesis on the $(-1)-$special system, we know that $\nu({\mathcal{M}})\geq 0$.
Moreover the  $C_i$'s are fixed for $i=s+1, \dots, t$, hence one has
$$\nu(\LL-\sum_{i=1}^sN_iC_i)=\nu(\LL-\sum_{i=1}^tN_iC_i)=\nu({\mathcal{M}})\geq 0$$
and we can conclude that $C=\sum_{i=1}^sN_iC_i$ is an $(N_1, \dots, N_s)-$special 
effect configuration for $\LL$. Then $\LL$ is numerically special.
\vskip0.3cm
\noindent
{\bf[(HHC) $\Rightarrow$ (CSEC)]} Suppose that the Harbourne--Hirschowitz Conjecture holds. 
As in the previous case, we prove that a $(-1)-$curve appearing in a $(-1)-$special 
system and splitting off with at least multiplicity two is an $h^1-$special effect curve.
Let $\LL$ be a special system. Then there is at least a $(-1)-$curve $C$ 
such that $\LL\cdot C<-N$, $N>1$. Then $h^0(\LL_{|C})=0$ and, by Riemann--Roch, $h^1(\LL_{|C})=N-1>0$ 
so that conditions $(a)$ and $(c)$ of Definition \ref{h1SEV} are satisfied. At this point it is 
important to observe that $\LL-C$ could be special. However the speciality of $\LL-C$ has no effect on $h^0(\LL-C)$. 
In fact if $\LL-C$ is non-special then $\LL-C$ contains the residual system ${\mathcal{M}}$ and, 
by definition of $(-1)-$special system, $\nu({\mathcal{M}})\geq 0$ so that $h^0(\LL-C)\not=0$. 
If $\LL-C$ is special, then, by (\ref{h0h1}) we surely have $h^0(\LL-C)\not=0$. 
Hence condition $(b)$ is satisfied.
\vskip0.3cm
\noindent
{\bf [(NSEC) $\Rightarrow$ (SC)]} 
Suppose that the Numerical Special Effect Conjecture holds. 
Let $\LL$ be a special system, then there is an $(\alpha_1, \dots, \alpha_t)-$
special effect configuration or an $\alpha-$special effect curve for $\LL$. 
We prove only the case in which there is a special effect 
configuration for $\LL$, being the other one similar.

Let $X=\sum_{i=1}^t\alpha_i Y_i$  be the special effect configuration. It is enough to fix 
our attention on $Y_1$. Since, by hypothesis,
 $\nu(\LL-\sum_{i=1}^t\alpha_i Y_i)\geq 0$, one has  
$h^0(\LL-\alpha_1 Y_1)\geq 1$ and 
we can apply Lemma \ref{inter}. Thus $Y_1^2\leq -1$ and, by Lemma \ref{LY}, we have
$$\LL \cdot Y_1 <  \frac{(\alpha_1+1)}{2}\tilde{Y_1}^2<-1.$$
Thus $Y_1$ is a fixed multiple component of $\LL$ and Segre's Conjecture holds.
\vskip0.3cm
\noindent
{\bf [(CSEC) $\Rightarrow$ (SC)]} 
Suppose that the Cohomological Special Effect Conjecture holds. 
Let $\LL$ be a special system, then there exists 
an $h^1-$special effect curve $Y$ for $\LL$. 
By condition $(b)$ of Definition \ref{h1SEV} we know that $Y$ splits from $\LL$,
then it is enough to show that $Y$ splits off 
at least with multiplicity $2$. 
Since $Y$ is irreducible, we have $Y^2 \geq g-1$ where $g$ is the genus of $Y$. 
By Riemann--Roch and $h^0(\LL_{|Y})=0$  we have $\LL\cdot Y = g-1-h^1(\LL_{|Y})$. 
Then we compute
$$(\LL-Y)\cdot Y = \LL \cdot Y - Y^2 \leq g-1-h^1(\LL_{|Y})-(g-1)=-h^1(\LL_{|Y})<0$$
and the claim follows.
\end{proof}


\section{First examples of special effect varieties in higher dimension}\label{highexamples}

Since a curve in $\PP^2$ is also a divisor, when we pass to analyze
the case of special linear systems in $\PP^n$, $n \geq
3$, we can pose the question if it is natural to consider special
effect varieties of every codimension (i.e. not only curves or not
only divisors). This more general situation is justified in
\cite{higherSEV}, where we prove, for example, that $\PP^s$, $1 \leq s \leq n-1$ can
be a special effect varieties for a given system $\LL$.

The definition of a special effect variety
$Y$ such that codim$(Y,\PP^n)\leq n-2$ is more difficult than the
codimension one case. Thus, here, we consider only when $Y$ is a
divisor. Obviously, in this situation, Definitions \ref{aSEP},  \ref{aSEV},
\ref{aSEC} and  \ref{h1SEV} remain the same.

\begin{Ex}\label{ex1}\rm Let $\LL$ be the system 
  $\LL_{3,4}(2^9)$ in Theorem \ref{P^nspecial}. Consider a quadric
  $Q\subset \PP^3$ through the nine points of $\LL$.
Obviously $\nu(|Q_3|)=0$
Moreover, one has 
$$\nu(\LL-Q)=\nu(\LL-2Q)=0$$
while 
$$\nu(\LL)=-2.$$
Thus $Q$ is a $2-$special effect variety (hypersurface) for  $\LL_{3,4}(2^9)$.
\end{Ex}

\begin{Ex}\rm In the same way we can prove that the quadric $Q \subset
  \PP^4$ 
is a $2-$special effect variety for  $\LL_{4,4}(2^{14})$.
\end{Ex}

\begin{Ex} Consider again the situation of Example \ref{ex1}. We prove
  that $Q$ is an $h^1-$special effect variety for $\LL$. Since
  $\LL(-Q)\cong \OO_{\PP^3}(Q)$ one has
$$H^0(\LL-Q)=1 \mbox{  and  } H^i(\LL-Q)=0, \quad i\geq 1$$
and condition $(b)$ is satisfied. Since we know that 
$$H^0(\LL)=1, \quad H^1(\LL)=2 \mbox{  and  } H^i(\LL)=0, \quad i\geq 2$$
we can conclude that
$$H^0(\LL_{|Q})=0 \mbox{  and  } H^1(\LL_{|Q})=2.$$
Thus conditions $(a)$ and $(c)$ hold and the claim follows.
\end{Ex}

\begin{Ex} In the same way we can prove that the quadric $Q \subset \PP^4$ 
is an $h^1-$special effect variety for  $\LL_{4,4}(2^{14})$.
\end{Ex}

\begin{Rmk}
In the previous examples we shows that the quadrics are both $\alpha-$
and $h^1-$special effect varieties for the same system. This is not
true in general. In fact, in \cite{higherSEV} we show that a plane
$\pi\subset \PP^3$ is
a $1-$special effect variety for $\LL:=\LL_{3,6}(4^3)$, but it is not an
$h^1-$special effect variety for the same system. 
\end{Rmk}


\section{Special effect curves on surfaces}\label{surfaces}

It could be interesting to extend the concept of special effect
curves to surfaces different from $\PP^2$.

We just give here some examples which show some important evidence.

\begin{Ex}{\bf Hirzebruch surfaces} \rm Let $\FF_e$, $e \ge 0$, be the Hirzebruch surface with invariant $e$, i.e. such that $-e$ is the minimal self-intersection of a section
of the ruling of $\FF_e$. 
 We have $\mbox{Pic}(\FF_e)\cong {\mathbb{Z}}\oplus {\mathbb{Z}}$ and
 we take, 
as a basis of $\mbox{Pic}(\FF_e)$, a section $h$ of the ruling
$f: \FF_e \to {\bf {P}}^1$ with $h^2=-e$ and a class, $F$, of
 $f$. Thus $h\cdot F = 1$ and $F^2 = 0$. 
The dimension of $H^0(\FF_e,\mathcal {O}_{\FF_e}(ah+bF))$ is given by
\[
\begin{cases}
0 & \text{if $a \ge 0$ and $b < 0$}, \\
\sum _{i=0}^{t-1}(b-ie+1) & \text{if $0 \le b < te$ for some $t\in
  {\mathbb{Z}}$,  with $0 \le t \le a$} \\
\frac{(2b+2-ae)(a+1)}{2} & \text{if $a \ge 0$ and $b \ge ae-1$}
\end{cases}
\]
and $h^1(\FF_e,\mathcal {O}_{\FF_e}(ah+bF)) = 0$ if $a \ge 0$ and $b \ge ae-1$. \\
We denote a system on $\FF_e$ by $\LL(a,b):=|ah+bF|$. 

Laface, in \cite{Laface}, gives a different definition of $(-1)-$special system. For that, we need the following procedure.
\vskip0.1cm
Given a linear system $\LL:=|ah+bF-\sum_{i=1}^{h}m_iP_i|$ on $\FF_e$
\begin{itemize}
\item[1)] if it does exist a $(-1)-$curve $E$ such that $-t:=\LL\cdot E <0$ then substitute $\LL$ with $\LL-tE$ and go to step 1), else go to step 2).
\item[2)] if $\LL\cdot h <0$ then substitute $\LL$ with $\LL-h$ and go to step 1), else finish.
\end{itemize}
After a finite number of steps, we have a new linear system ${\mathcal{M}}$, i.e., the residual linear system.

\begin{Def}\label{Lafacedef}
Let $\LL:=|ah+bF-\sum_{i=1}^{h}m_iP_i|$ and ${\mathcal{M}}$ on $\FF_e$ as above. Then $\LL$ is $(-1)-${\bf special} if $v({\mathcal{M}})>v(\LL)$.
\end{Def}

Then we can state again a modified Harbourne--Hirschowitz Conjecture:

\begin{Conj}[A. Laface, 2002]\label{LafConj}
A system $\LL(-\sum_{i=1}^hm_iP_i)$ on a ${\mathbb{F}}_e$  is special if and only if is $(-1)-$special.
\end{Conj}

This time, for the speciality of a linear system $\LL$ such that
$\LL=\sum_{i=1}^tN_iC_i + {\mathcal{M}}$ it is not enough to have
$v({\mathcal{M}})\ge 0$ and $N_i\ge2$ for at least one index $i$. 
Following the argument of the main theorem in \cite{Laface} It is easy
to construct several examples of special system in $\FF_e$, $e \geq
4$,
such that the Harbourne--Hirschowitz does not hold (see Example 3.4.4 in \cite{Bocci}).

The interested reader can look at Laface's article for a deep understanding. We just recall the main results contained in it.

\begin{Prop}\label{ListaLaface}
Denote by  $\LL_e(a,b,m^h)$ the system $\LL_e(a,b)(-\sum_{i=1}^hmP_i)$. All homogeneous $(-1)-$special systems with multiplicity $m \leq 3$ on $\FF_e$ are listed in the following table:
\begin{center}
\begin{tabular}{ccc}
\hline
system & virtdim$(\LL)$ & dim$(\LL)$ \\\hline
$\LL_1(4,4,2^5)$ & $-1$ & $0$\\
$\LL_1(6,6,3^5)$ & $-3$ & $0$\\
$\LL_5(4,21,3^{10})$ & $-1$ & $0$\\
$\LL_6(4,24,3^{11})$ & $-1$ & $0$\\
$\LL_e(2,2d+2e,2^{2d+e+1})$ & $-1$ & $0$\\
$\LL_e(0,d,2^r)$ & $d-3r$ & $d-2r$\\
$\LL_e(2,4d+3e+1,3^{2d+e+1})$ & $-1$ & $0$\\
$\LL_e(3,3d+3e+1,3^{2d+e+1})$ & $1$ & $2$\\
$\LL_e(3,3d+3e,3^{2d+e+1})$ & $-3$ & $0$\\
$\LL_e(1,d+e,3^r)$ & $2d+e-6r+1$ & $2d+e-5r+1$\\
$\LL_e(0,d,3^r)$ & $d-6r$ & $d-3r$\\
\hline
\end{tabular}
\end{center}
\end{Prop}

\begin{Thm}
Every special homogeneous system of multiplicity $\leq 3$ on a $\FF_e$ surface is a $(-1)-$special system.
\end{Thm}

After we modify the condition for $\alpha$ by respect to the degree of
$\LL$ and $Y$, we can give again the definition of $\alpha-$special
effect property and arrive again to state the Numerical Special Effect
Conjecture. One has the following

\begin{Thm}
The Numerical Special Effect
Conjecture on Hirzebruch surface holds for all special systems listed
in Proposition \ref{ListaLaface}.
\end{Thm}

\begin{proof}
It is enough to check by hand every single case on the previous
table. As an example we prove the case $\LL:=\LL_e(0,d,2^r)$, $d\geq 2r$.
Consider the curve $Y_1$ of bidegree $(0,1)$ passing through one of the
$r$ points in $\LL$, i.e. $Y_1$ corresponds to the system $\LL_e(0,1,1)$
Thue one has $\nu(\LL)=d-3r$ and $\nu(\LL-Y_1)=\nu(\LL-2Y_1)=d-3r+1$.
If $d-3r+1 \geq 0$ we conclude that $Y_1$ is a $2-$special effect curve
for $\LL$. In the other case we pass to study the system
$\LL':=\LL_e(0,d-2,2^{r-1})$ and we consider a new curve $Y_2$ passing
through one of the $r-1$ points of $\LL'$. As in the case of $Y_1$ we
conclude that $Y_2$ is a $2-$special effect curve
for $\LL'=\LL-2Y_1$. Going furhter we will obtain a $(2, \dots, 2)-$special
effect configuration $X=\sum_{i=1}^r2Y_i$ for $\LL$.
\end{proof}

Consider now the Cohomological Special Effect Conjecture.
Unluckily it does not hold for all special systems listed
in Proposition \ref{ListaLaface}. 

In fact, let $\LL$ be the special system $\LL_6(4,24,3^{11})$. We
know, by \cite{Laface},
that $\LL$ splits as $3E+h$, where $E$ is the $(-1)-$curve
corresponding to the system $\LL_6(1,8,1^{11})$.
By condition $h^0(\LL_{|Y})=0$, we know that an $h^1-$special effect
variety must split from $\LL$. Thus only $E$ and $h$ are the candidate
to be $h^1-$special for $\LL$.
Since $\LL\cdot E=-1$ (in fact $h$ ``hides'' the effective
multiplicity of $E$, see \cite{Bocci} or \cite{Laface}) one has
$h^1(\LL_{|E})=0$. Similarly, since $\LL\cdot h=0$, we have again
$h^1(\LL_{|h})=0$. 
Thus condition $(c)$ is never satisfied and both $E$ and $h$ are not $h^1-$special effect curves for $\LL$.
\end{Ex}

\begin{Ex}{\bf (K3 surfaces)} \rm 
Let $X$ be a K3 surface with $n=H^2\in2{\mathbb{Z}}$. Let $\LL:=\LL^n(d,m_1,\dots, m_h)$ be the system of curves $|dH|$ passing through points $P_1, \dots , P_h$ in general position on $X$ with multiplicities at least $m_1, \dots, m_h$. The virtual dimension of $\LL$ is given by
$$\nu(\LL)=d^2\frac{H^2}{2}-\sum_{i=1}^h\frac{m_i(m_i+1)}{2}+1.$$
In \cite{DeVolLaf},  De Volder and  Laface state a conjecture for linear systems on a K3 surface and, moreover, they proved it is equivalent to the Segre Conjecture. i.e. if $\LL$ is special on $X$ then $\LL$ has a multiple fixed component.

\begin{Conj}[De Volder--Laface]\label{DVLC} Let $\LL$ and $X$ be as above.
\begin{itemize}
\item[(i)] $\LL$ is special if and only if $\LL=\LL^4(d,2d)$ or $\LL=\LL^2(d,d^2)$ with $d\geq 2$;
\item[(ii)] if $\LL$ is non-empty then its general divisor has exactly the imposed multiplicities at the points $P_i$;
\item[(iii)] if $\LL$ is non-special and has a fixed irreducible component $C$ then
\begin{itemize}
\item[a)] $\LL:=\LL^2(m+1,m+1,m)=mC+\LL^2(1,1)$ with $C=\LL^2(1,1^2)$ or
\item[b)] $\LL=2C$, $C\in \{\LL^4(1,1^3),\LL^6(1,1,2),\LL^{10}(1,3) \}$ or
\item[c)] $\LL=C$. 
\end{itemize}
\item[(iv)] if $\LL$ has no fixed component then either its general element is irreducible or $\LL=\LL^2(2,2)$.
\end{itemize}
\end{Conj}
\vskip0.2cm
Consider the system $\LL=\LL^2(d,d^2)$. Its virtual dimensions is
$$\nu(\LL)=d^2-d(d+1)+1=1-d.$$
Let $C_1$ be the curve $\LL^2(1,1^2)$, then $C_1$ is a $d-$special
effect curve for $\LL$ since
$\nu(\LL-dC_1)=0$.
In a similar way we can prove that $C_2:=\LL^4(1,2)$
is a $d-$special
effect curve  $\LL^4(d,2d)$. Moreover, we can see that 
$\nu(\LL-C)=\nu(\LL)$ when $C$ is one of the curve in cases
$(iii)a)-c)$ of the conjecture and $\LL$ is the relative system to $C$. 
\vskip0.2cm
Passing to the $h^1-$special effect curves, we can observe
that that $C_1:=\LL^2(1,1^2)$ and $C_2:=\LL^4(1,2)$ are genus
two curves with self-intersection equal to zero.
Applying Riemann--Roch we discover that $h^0(\LL_{|C_t})=0$ and
$h^1(\LL_{|C_t})=1$, where $\LL$ is the relative system to $C_t$ in case $(i)$ in
Conjecture \ref{DVLC} ($t=1,2$). Since $h^0(\LL-C_i)>0$ we conclude that systems
in $(i)$ are cohomologically special.
Finally we can see that no curve $C$ in cases $(iii)a) - c)$ are
$h^1-$special effect curve. In fact, in all case in $(iii)a)-b)$ one
has $h^i(\LL_{|C})=0$, $i=0,1$. While the cuve in $(iii)c)$ does not
fit the hypothesis in Definition \ref{h1SEV}, since $\LL\cong \OO_X(C)$.

Thus we can state the following
\begin{Thm} Conjecture \ref{DVLC}
implies both Numerical and Cohomological Conjectures.
\end{Thm}
\end{Ex}

\end{document}